\documentclass[12pt]{article}

\renewcommand{\leq}{\leqslant}
\renewcommand{\geq}{\geqslant}

\usepackage{amsthm,amsmath,amssymb,latexsym}

\usepackage{epsfig,cite,lineno}

\newtheorem{thm}{Theorem}[section]
\newtheorem{prop}[thm]{Proposition}

\newtheorem{cor}[thm]{Corollary}

\newtheorem{lem}[thm]{Lemma}
\newtheorem{conj}[thm]{Conjecture}

\theoremstyle{remark}

\def\Z{\mathbb{Z}^{+}}
\def\N{\mathbb{N}}

\def\pf{\noindent {\it Proof.} }

\numberwithin{equation}{section}
\renewcommand{\qed}{\hfill$\Box$\medskip}

\begin{document}

\begin{center}
{\Large\bf Factors of sums and alternating sums of products of\\[5pt]
$q$-binomial coefficients and powers of $q$-integers}
\end{center}
\vskip 2mm \centerline{Victor J. W. Guo\footnote{Corresponding author} and Su-Dan Wang$^{2}$}
\begin{center}
{\footnotesize $^1$School of Mathematical Sciences, Huaiyin Normal University, Huai'an, Jiangsu 223300,\\
 People's Republic of China\\
{\tt jwguo@hytc.edu.cn}\\[10pt]
%$^2$Universit\'e de Lyon, Lyon, F-69003, France\\

$^2$Department of Mathematics, East China Normal University, Shanghai 200062,\\
 People's Republic of China\\
 {\tt sudan199219@126.com}
 }
\end{center}

%\centerline{September 3, 2009}

\vskip 0.7cm {\noindent{\bf Abstract.}
We prove that, for all positive integers
$n_1, \ldots, n_m$, $n_{m+1}=n_1$, and non-negative integers $j$ and $r$ with $j\leqslant m$, the following two expressions
\begin{align*}
&\frac{1}{[n_1+n_m+1]}{n_1+n_{m}\brack n_1}^{-1}\sum_{k=0}^{n_1} q^{j(k^2+k)-(2r+1)k}[2k+1]^{2r+1}\prod_{i=1}^m {n_i+n_{i+1}+1\brack n_i-k},\\[5pt]
&\frac{1}{[n_1+n_m+1]}{n_1+n_{m}\brack n_1}^{-1}\sum_{k=0}^{n_1}(-1)^k q^{{k\choose 2}+j(k^2+k)-2rk}[2k+1]^{2r+1}\prod_{i=1}^m {n_i+n_{i+1}+1\brack n_i-k}
\end{align*}
are Laurent polynomials in $q$ with integer coefficients, where $[n]=1+q+\cdots+q^{n-1}$ and ${n\brack k}=\prod_{i=1}^k(1-q^{n-i+1})/(1-q^i)$.
This gives a $q$-analogue of some divisibility results of sums and alternating sums involving binomial coefficients and powers of integers
obtained by Guo and Zeng. We also confirm some related conjectures of Guo and Zeng by establishing their $q$-analogues.
Several conjectural congruences for sums involving products of $q$-ballot numbers $\left({2n\brack n-k}-{2n\brack n-k-1}\right)$
are proposed in the last section of this paper.
}



%%%%%%%%%%%%%%%%%%%%%%%%%%%%%%%%%%%%%%%%%%%%%%%%%%%%%%%%%%%%%%%%%%%%%%%%%%%%%%%%%%%%%%%%%%%%%%%%%%%
\section{Introduction}
In 2011, the first author and Zeng \cite{GZ2011} prove that, for all positive integers
$n_1, \ldots, n_m$, $n_{m+1}=n_1$, and any non-negative integer $r$, there holds
\begin{align}
\sum_{k=0}^{n_1}\varepsilon^k (2k+1)^{2r+1}\prod_{i=1}^{m} {n_i+n_{i+1}+1\choose n_i-k} \equiv 0 \mod (n_1+n_m+1){n_1+n_m\choose n_1}, \label{eq:guozeng-1}
\end{align}
where $\varepsilon=\pm 1$. The congruence \eqref{eq:guozeng-1} is very similar to the following congruences:
\begin{align}
\sum_{k=-n_1}^{n_1}(-1)^k\prod_{i=1}^m
{n_i+n_{i+1}\choose n_i+k}\equiv 0\mod{ {n_1+n_m\choose n_1} }, \label{eq:guozeng-2} \\[5pt]
2\sum_{k=1}^{n_1}k^{2r+1}\prod_{i=1}^{m} {n_i+n_{i+1}\choose n_i+k}\equiv 0\mod{ n_1{n_1+n_m\choose n_1} },  \label{eq:guozeng-3}
\end{align}
where $n_{m+1}=n_1$, which were obtained by Guo, Jouhet, and Zeng \cite{GJZ}, and Guo and Zeng \cite{GZ2010}, respectively.
Note that \eqref{eq:guozeng-2} is a generalization of the following congruence due to Calkin \cite{Calkin}:
\begin{align*}
\sum_{k=-n}^n (-1)^k{2n\choose n+k}^m\equiv 0\mod{{2n\choose n}}\qquad \textrm{for $m\geq 1$.}  %\label{eq:Calkin}
\end{align*}
It is known that both \eqref{eq:guozeng-2} and \eqref{eq:guozeng-3} have neat $q$-analogues (see \cite{GJZ} and \cite{GW}).
It is also worth mentioning that $q$-analogues of classical congruences have been widely studied during the last decade (see, for example,
\cite{PS,SP,Tauraso2012,Tauraso2013}).

The first aim of this paper is to give a $q$-analogue of \eqref{eq:guozeng-1}. Recall that
the {\it $q$-integers} are defined as $[n]=1+q+\cdots+q^{n-1}$ and the {\it $q$-binomial coefficients} are defined by
$$
{n\brack k}=
\begin{cases}\displaystyle\prod_{i=1}^k\frac{1-q^{n-i+1}}{1-q^i} &\text{if $k\geqslant 0$,} \\[10pt]
0 &\text{otherwise.}
\end{cases}
$$
Let $D$ be a polynomial in $q$. We say that two Laurent polynomials $A$ and $B$ in $q$ are congruent modulo $D$, denoted by
$A\equiv B\mod D$, if $(A-B)/D$ is still a Laurent polynomial in $q$. Let $\N$ denote the set of non-negative integers and
$\Z$ the set of positive integers. Our first result is as follows.
\begin{thm}\label{thm:fatorodd-first}
Let $n_1,\ldots,n_{m}\in\Z$, $n_{m+1}=n_1$, and $j,r\in\N$ with $j\leqslant m$. Then modulo
$[n_1+n_m+1]{n_1+n_{m}\brack n_1}$,
\begin{align}
&\sum_{k=0}^{n_1} q^{j(k^2+k)-(2r+1)k}[2k+1]^{2r+1}\prod_{i=1}^m {n_i+n_{i+1}+1\brack n_i-k}\equiv 0, \label{eq:first-1} \\[5pt]
&\sum_{k=0}^{n_1}(-1)^k q^{{k\choose 2}+j(k^2+k)-2rk}[2k+1]^{2r+1}\prod_{i=1}^m {n_i+n_{i+1}+1\brack n_i-k} \equiv 0.  \label{eq:first-2}
\end{align}
\end{thm}

The first author and Zeng \cite{GZ2011} also prove that, for all positive integers
$n_1, \ldots, n_m$, $n_{m+1}=n_1$, and any non-negative integer $r$,
\begin{align}
&\hskip -3mm \sum_{k=0}^{n_1}k^r(k+1)^r(2k+1)\prod_{i=1}^{m} {n_i+n_{i+1}+1\choose n_i-k} \notag\\[5pt]
&\equiv 0 \mod (n_1+n_m+1){n_1+n_m\choose n_1}n_1^{\min\{1,r\}}n_m^{\min\{1,{r\choose 2}\}}, \label{eq:guozeng-4}\\[5pt]
&\hskip -3mm \sum_{k=0}^{n_1}(-1)^{k} k^r(k+1)^r(2k+1)\prod_{i=1}^{m} {n_i+n_{i+1}+1\choose n_i-k} \notag\\[5pt]
&\equiv 0 \mod (n_1+n_m+1){n_1+n_m\choose n_1}n_1 ^{\min\{1,r\}}n_m ^{\min\{1,r\}}.   \label{eq:guozeng-5}
\end{align}
Actually in \cite{GZ2011} the congruence \eqref{eq:guozeng-1} is deduced from \eqref{eq:guozeng-4} and \eqref{eq:guozeng-5} by noticing that
$$
(2k+1)^{2r}=(4k^2+4k+1)^r =\sum_{i=0}^r{r\choose i}4^{i}k^i(k+1)^i.
$$

The second aim of this paper is to give the following $q$-analogue of \eqref{eq:guozeng-4} and \eqref{eq:guozeng-5}.

\begin{thm}\label{thm:fatorodd-second}
Let $n_1,\ldots,n_{m}\in\Z$, $n_{m+1}=n_1$, and $j,r\in\N$ with $j\leqslant m$. Then
\begin{align*}
&\hskip -2mm \sum_{k=0}^{n_1} q^{j(k^2+k)-(r+1)k}[2k+1][k]^r[k+1]^r\prod_{i=1}^m {n_i+n_{i+1}+1\brack n_i-k} \notag\\[5pt]
&\equiv 0 \mod [n_1+n_m+1]{n_1+n_m\brack n_1}[n_1]^{\min\{1,r\}}[n_m]^{\min\{1,{r\choose 2}\}}, \\[5pt]
&\hskip -2mm \sum_{k=0}^{n_1}(-1)^k q^{{k\choose 2}+j(k^2+k)-rk} [2k+1][k]^r[k+1]^r\prod_{i=1}^m {n_i+n_{i+1}+1\brack n_i-k} \notag \\[5pt]
&\equiv 0 \mod [n_1+n_m+1]{n_1+n_m\brack n_1}[n_1]^{\min\{1,r\}}[n_m]^{\min\{1,r\}}.
\end{align*}
\end{thm}

Not like the $q=1$ case, it seems that Theorem \ref{thm:fatorodd-first} cannot be derived from Theorem \ref{thm:fatorodd-second} directly.

The $q$-ballot numbers $A_{n,k}(q)$ ($0\leq k\leq n$) are defined by
\begin{align}
A_{n,k}(q)=q^{n-k}\frac{[2k+1]}{[2n+1]}{2n+1\brack n-k}={2n\brack n-k}-{2n\brack n-k-1}.
\end{align}
Note that sums involving the ballot numbers $A_{n,k}:=A_{n,k}(1)$ have been considered by Miana and Romero \cite[Theorem 10]{MR2}
and Guo and Zeng \cite{GZ2011}.

The third aim of this paper is to give the following congruences involving $q$-ballot numbers. Note that the $q=1$ case
confirms a conjecture of Guo and Zeng \cite[Conjecture 1.3]{GZ2011}.
\begin{thm}\label{thm:q-ballot}
Let $n,s\in\Z$ and $r,j\in\N$ with $r+s\equiv 1\pmod{2}$ and $j\leqslant s$. Then
\begin{align}
&\sum_{k=0}^{n}q^{j(k^2+k)-rk}[2k+1]^r A_{n,k}(q)^s\equiv 0 \mod{ {2n\brack n} }, \label{eq:q-ballt-1}\\[5pt]
&\sum_{k=0}^{n}(-1)^k q^{{k\choose 2}+j(k^2+k)-(r-1)k}[2k+1]^r A_{n,k}(q)^s  \equiv 0 \mod{ {2n\brack n}}.  \label{eq:q-ballt-2}
\end{align}
\end{thm}

Let $[n]!=[n][n-1]\cdots[1]$ be the {\it $q$-factorial} of $[n]$. It is easy to see that, for all $m,n\in\N$, the expression
$\frac{[2m]![2n]!}{[m+n]![m]![n]!}$ is a polynomial in $q$ by writing a $q$-factorial as
a product of cyclotomic polynomials. The polynomials $\frac{[2m]![2n]!}{[m+n]![m]![n]!}$ are usually called the {\it $q$-super Catalan numbers}.
Warnaar and Zudilin \cite[Proposition 2]{WZ} have shown that the $q$-super Catalan numbers are polynomials in $q$
with non-negative integer coefficients.

We shall also prove the following congruences modulo $q$-super Catalan numbers.
\begin{thm}\label{thm:mn-qq}
Let $m,n,s,t\in\Z$ and $j,r\in \N$ with $r+s+t\equiv 1\pmod 2$ and $j\leqslant s+t$. Then
\begin{align*}
[m+n+1]\sum_{k=0}^m q^{j(k^2+k)-rk} [2k+1]^{r}A_{m,k}(q)^s  A_{n,k}(q)^t
&\equiv 0 \mod \frac{[2m]![2n]!}{[m+n]![m]![n]!}, \\[5pt]
[m+n+1]\sum_{k=0}^m (-1)^kq^{{k\choose 2}+j(k^2+k)-(r-1)k} [2k+1]^{r}A_{m,k}(q)^s  A_{n,k}(q)^t
&\equiv 0 \mod \frac{[2m]![2n]!}{[m+n]![m]![n]!}.
\end{align*}
\end{thm}

Note that the $q=1$ case of Theorem \ref{thm:mn-qq} confirms another conjecture of Guo and Zeng \cite[Conjecture 6.10]{GZ2011}.
It should also be mentioned that Theorem \ref{thm:mn-qq} in the case where $m=n$ gives the $s\geqslant 2$ case
of Theorem \ref{thm:q-ballot} (see \eqref{eq:gcd-bino}).

The paper is organized as follows. We shall prove Theorem \ref{thm:fatorodd-first} for $m=1$ in Section 2 and prove Theorem \ref{thm:fatorodd-second} for $m=1$
in Section 3. A proof of Theorems \ref{thm:fatorodd-first} and \ref{thm:fatorodd-second} for $m\geqslant 2$ will be given in Section 4.
The $q$-Chu-Vandermonde identity and the $q$-Dixon identity will play a key role in our proof. We shall prove Theorems \ref{thm:q-ballot} and \ref{thm:mn-qq} in Sections 5
and 6, respectively. We give some consequences of Theorem \ref{thm:fatorodd-first} and some related conjectures in Section~7.

%%%%%%%%%%%%%%%%%%%%%%%%%%%%%%%%%%%%%%%%%%%%%%%%%%%%%%%%%%%%%%%%%%%%%%%%%%%%%%%%%%%%%%%%%%%%%%%%%%%%%%%%%%%%%%%%%%%%%
\section{Proof of Theorem \ref{thm:fatorodd-first} for $m=1$}
The {\it $q$-shifted factorials} (see \cite{GR}) are defined as $(a;q)_0=1$
and $(a;q)_n=(1-a)(1-aq)\cdots (1-aq^{n-1})$ for $n=1,2,\ldots.$
In order to prove Theorem \ref{thm:fatorodd-first} for $m=1$, we shall first establish the following result.
\begin{lem}Let $n\in\Z$ and $s\in\N$. Then
\begin{align}
&\sum_{k=0}^{n}q^{-k}[2k+1] {2n+1\brack n-k}(q^{-k};q)_s (q^{k+1};q)_s =(-1)^s q^{{s\choose 2}-sn-n}[2n+1]{2n\brack n}{n\brack s} (q;q)_s^2, \label{eq:r1m1-1}\\[5pt]
&\sum_{k=0}^{n}q^{k^2}[2k+1] {2n+1\brack n-k}(q^{-k};q)_s (q^{k+1};q)_s =(-1)^s q^{s\choose 2}[2n+1]{2n\brack n}{n\brack s} (q;q)_s^2, \label{eq:r1m1-2} \\[5pt]
&\sum_{k=0}^{n}(-1)^k q^{k\choose 2}[2k+1] {2n+1\brack n-k}(q^{-k};q)_s (q^{k+1};q)_s =0, \label{eq:r1m1-3}\\[5pt]
&\sum_{k=0}^{n}(-1)^k q^{\frac{3k^2+k}{2}}[2k+1] {2n+1\brack n-k}(q^{-k};q)_s (q^{k+1};q)_s
=q^{s^2} [2n+1]{2n\brack n}{n\brack s}(q;q)_n(q;q)_s.   \label{eq:r1m1-4}
\end{align}
\end{lem}
\pf We proceed by induction on $s$. For $s=0$, we have
\begin{align}
\sum_{k=0}^{n}q^{-k}[2k+1] {2n+1\brack n-k}
&=q^{-n}[2n+1]\sum_{k=0}^{n}\left({2n\brack n-k}-{2n\brack n-k-1}\right) \notag\\[5pt]
&=q^{-n}[2n+1]{2n\brack n}, \notag\\[5pt]
\sum_{k=0}^{n}q^{k^2}[2k+1] {2n+1\brack n-k}
&=[2n+1]\sum_{k=0}^{n}\left(q^{k^2}{2n\brack n-k}-q^{(k+1)^2}{2n\brack n-k-1}\right) \notag\\[5pt]
&=[2n+1]{2n\brack n}, \notag \\[5pt]
\sum_{k=0}^{n}(-1)^k q^{k\choose 2}[2k+1] {2n+1\brack n-k}
&=q^{-n}[2n+1]\sum_{k=0}^{n}(-1)^k q^{k+1\choose 2}\left({2n\brack n-k}-{2n\brack n-k-1}\right) \notag\\[5pt]
&=q^{-n}[2n+1]\sum_{k=-n}^{n}(-1)^k q^{k+1\choose 2}{2n\brack n-k} \notag\\[5pt]
&=0, \label{eq:app-qbino}
\end{align}
and
\begin{align}
\sum_{k=0}^{n}(-1)^k q^{\frac{3k^2+k}{2}}[2k+1] {2n+1\brack n-k}
&=[2n+1]\sum_{k=0}^{n}(-1)^k q^{k+1\choose 2}  \notag\\[5pt]
&\quad\times\left(q^{k^2}{2n\brack n-k}-q^{(k+1)^2}{2n\brack n-k-1}\right)\notag\\[5pt]
&=[2n+1]\sum_{k=-n}^{n}(-1)^k q^{\frac{3k^2+k}{2}}{2n\brack n-k}  \notag\\[5pt]
&=[2n+1]{2n\brack n}(q;q)_n,  \label{eq:app-qdixon}
\end{align}
where the equality \eqref{eq:app-qbino} follows from the $q$-binomial theorem (see \cite[p.~36, Theorem 3.3]{Andrews}):
\begin{align*}
(x;q)_{N} &=\sum_{k=0}^{N}(-1)^k q^{k\choose 2}{N\brack k} x^k
\end{align*}
by taking $x=q^{-n}$ and $N=2n$, while the equality
\eqref{eq:app-qdixon} is the $l,m\to\infty$ case of the $q$-Dixon identity:
\begin{align*}
\sum_{k=-n}^{n}(-1)^k q^{\frac{3k^2+k}{2}}{l+m\brack l+k}{m+n\brack m+k}{n+l\brack n+k}
=\frac{(q;q)_{l+m+n}}{(q;q)_l (q;q)_m (q;q)_n}
\end{align*}
(see \cite{GZ-Dixon} for a short proof).

Suppose that the identities \eqref{eq:r1m1-1}--\eqref{eq:r1m1-4} are true for $s$. Noticing the relation
\begin{align*}
&\hskip -2mm {2n+1\brack n-k}(q^{-k};q)_{s+1} (q^{k+1};q)_{s+1}  \\
&=(1-q^{s-n})(1-q^{s+n+1}){2n+1\brack n-k}(q^{-k};q)_{s} (q^{k+1};q)_{s}  \\
&\quad{}+q^{s-n}(1-q^{2n})(1-q^{2n+1}){2n-1\brack n-k-1}(q^{-k};q)_{s} (q^{k+1};q)_{s},
\end{align*}
we can easily deduce that the identities \eqref{eq:r1m1-1}--\eqref{eq:r1m1-4} hold for $s+1$.
\qed

\noindent{\it Remark.} We have the following generalization of \eqref{eq:r1m1-3}:
\begin{align*}
&\sum_{k=0}^{n}(-1)^k q^{k\choose 2}[2k+1] {2n+1\brack n-k}(xq^{-k};q)_s (xq^{k+1};q)_s \notag\\
&\quad{} =x^n q^{-n}[2n+1]{2n\brack n}{s\brack n}\frac{(x;q)_{s-n}(x;q)_{s+1}(q;q)_{n}^2}{(x;q)_{n+1}},
\end{align*}
which can be proved in the same way as before.

We shall prove Theorem {\rm\ref{thm:fatorodd-first}} for $m=1$ in the following more general form:
\begin{thm}
Let $n\in\Z$ and $r,s\in\N$. Then modulo $[2n+1]{2n\brack n}$,
\begin{align}
&\sum_{k=0}^{n}q^{-(2r+1)k}[2k+1]^{2r+1} {2n+1\brack n-k}(q^{-k};q)_s (q^{k+1};q)_s \equiv 0, \label{eq:mod-r1m1-1}\\[5pt]
&\sum_{k=0}^{n}q^{k^2-2rk}[2k+1]^{2r+1} {2n+1\brack n-k}(q^{-k};q)_s (q^{k+1};q)_s \equiv 0, \label{eq:mod-r1m1-2} \\[5pt]
&\sum_{k=0}^{n}(-1)^k q^{{k\choose 2}-2rk}[2k+1]^{2r+1} {2n+1\brack n-k}(q^{-k};q)_s (q^{k+1};q)_s \equiv 0, \label{eq:mod-r1m1-3}\\[5pt]
&\sum_{k=0}^{n}(-1)^k q^{\frac{3k^2+k}{2}-2rk}[2k+1]^{2r+1} {2n+1\brack n-k}(q^{-k};q)_s (q^{k+1};q)_s
\equiv 0.   \label{eq:mod-r1m1-4}
\end{align}
\end{thm}
\pf We proceed by induction on $r$.
Denote the left-hand side of \eqref{eq:mod-r1m1-1} by $A_{r}(n,s)$.
By \eqref{eq:r1m1-1}, we know that \eqref{eq:mod-r1m1-1} is true for $r=0$.
For $r\geqslant 1$, suppose that
$$
A_{r-1}(n,s)\equiv 0\mod [2n+1]{2n\brack n}
$$
holds for all non-negative integers $n$ and $s$.
It is easy to check that
\begin{align*}
{2n+1\brack n-k}[2k+1]^2
&=q^{2k-2n}{2n+1\brack n-k}[2n+1]^2 \\[5pt]
&\quad{}-q^{2k-2n}{2n-1\brack n-k-1}[2n][2n+1](1+q^{n-s})(1+q^{n+s+1}) \\[5pt]
&\quad{}+q^{2k-n-s}{2n-1\brack n-k-1}[2n][2n+1](1-q^{s-k})(1-q^{s+k+1}),
\end{align*}
and therefore,
\begin{align}
A_{r}(n,s)
&=q^{-2n}[2n+1]^2A_{r-1}(n,s)-q^{-2n}[2n][2n+1](1+q^{n-s})(1+q^{n+s+1})A_{r-1}(n-1,s) \notag\\[5pt]
&\quad{}+q^{-n-s}[2n][2n+1]A_{r-1}(n-1,s+1).  \label{eq:rec-3-term}
\end{align}

By the induction hypothesis, we have
\begin{align*}
[2n][2n+1]A_{r-1}(n-1,s)
&\equiv [2n][2n+1]A_{r-1}(n-1,s+1)  \\[5pt]
&\equiv 0 \mod [2n][2n+1][2n-1]{2n-2\brack n-1}.
\end{align*}
Noticing that $[2n][2n+1][2n-1]{2n-2\brack n-1}=[2n+1]{2n\brack n}[n]^2$, the recurrence \eqref{eq:rec-3-term}
immediately implies that \eqref{eq:mod-r1m1-1} holds for $r$. Similarly, we can prove \eqref{eq:mod-r1m1-2}--\eqref{eq:mod-r1m1-4}.
\qed

%%%%%%%%%%%%%%%%%%%%%%%%%%%%%%%%%%%%%%%%%%%%%%%%%%%%%%%%%%%%%%%%%%%%%%%%%%%%%%%%%%%%%%%%%%%%%%%%%%%%%%%%%%%%%%%%%%%%%
\section{Proof of Theorem \ref{thm:fatorodd-second} for $m=1$}
For convenience, let
\begin{align*}
&P_r(n,j):=\sum_{k=0}^{n} q^{j(k^2+k)-(r+1)k}[2k+1][k]^r[k+1]^r {2n+1\brack n-k},\\[5pt]
&Q_r(n,j):=\sum_{k=0}^{n}(-1)^k q^{{k\choose 2}+j(k^2+k)-rk}[2k+1][k]^r[k+1]^r {2n+1\brack n-k}.
\end{align*}
Then the $m=1$ case of Theorem \ref{thm:fatorodd-second} can be restated as follows.
\begin{thm}\label{thm:fatorodd-secondm=1}
Let $n\in\Z$ and $r\in\N$. Then for $j=0,1$, there hold
\begin{align}
&P_r(n,j)\equiv 0 \mod{[2n+1]{2n\brack n}[n]^{\min\{2,r\}}},  \label{eq:srnj-0}\\[5pt]
&Q_r(n,j)\equiv 0 \mod{[2n+1]{2n\brack n}[n]^{\min\{2,2r\}}}.  \label{eq:trnj-0}
\end{align}
\end{thm}
\pf We proceed by induction on $r$. For $r=0$, by \eqref{eq:r1m1-1}--\eqref{eq:r1m1-4}, we have
\begin{align*}
&P_0(n,0)=q^{-n}[2n+1]{2n\brack n},\quad  P_0(n,1)=[2n+1]{2n\brack n},\\
&Q_0(n,0)=0\ (n\geqslant 1),\quad Q_0(n,1)=[2n+1]{2n\brack n}(q;q)_n.
\end{align*}
For $r\geq1$, observing that
\begin{align*}
q^{n-k}[k][k+1]{2n+1\brack n-k}=[n][n+1]{2n+1\brack n-k}-[2n][2n+1]{2n-1\brack n-k-1},
\end{align*}
we have the following recurrences:
\begin{align}
&P_r(n,j)=q^{-n}[n][n+1]P_{r-1}(n,j)-q^{-n}[2n][2n+1]P_{r-1}(n-1,j), \label{eq:rec-srnj}\\[5pt]
&Q_r(n,j)=q^{-n}[n][n+1]Q_{r-1}(n,j)-q^{-n}[2n][2n+1]Q_{r-1}(n-1,j)  \label{eq:rec-trnj}
\end{align}
for $n\geq1$.  From \eqref{eq:rec-srnj}--\eqref{eq:rec-trnj} we immediately get
\begin{align*}
&P_1(n,0)=q^{-2n}[n][2n+1]{2n\brack n},\quad P_2(n,0)=q^{-3n}[2][n]^2[2n+1]{2n\brack n},\\[5pt]
&P_1(n,1)=[n][2n+1]{2n\brack n},\quad P_2(n,1)=q^{-1}[2][n]^2[2n+1]{2n\brack n},\\[5pt]
&Q_1(1,0)=-q^{-1}[2][3],\quad Q_1(n,0)=0\ (n\geqslant 2),\quad Q_1(n,1)=-q[2n+1]{2n\brack n}[n]^2(q;q)_{n-1}.
\end{align*}
Therefore, the congruence \eqref{eq:srnj-0} is true for $r=0,1,2$, while the congruence \eqref{eq:trnj-0} is true for $r=0,1$.
We now assume that $r\geq3$ and \eqref{eq:srnj-0} holds for $r-1$ and $j=0,1$. Namely,
\begin{align*}
P_{r-1}(n,j)\equiv 0\mod[2n+1]{2n\brack n}[n]^2.
\end{align*}
It follows that
\begin{align*}
[2n][2n+1]P_{r-1}(n-1,j) \equiv 0\mod [2n][2n+1][2n-1]{2n-2\brack n-1}[n-1]^2.
\end{align*}
Since $[2n][2n+1][2n-1]{2n-2\brack n-1}=[2n+1]{2n\brack n}[n]^2$, from \eqref{eq:rec-srnj} we deduce that
$$
P_r(n,j)\equiv 0 \pmod{[2n+1]{2n\brack n}[n]^{2}}.$$
This completes the inductive step of \eqref{eq:srnj-0}. The proof of \eqref{eq:trnj-0} is exactly the same.
\qed

%%%%%%%%%%%%%%%%%%%%%%%%%%%%%%%%%%%%%%%%%%%%%%%%%%%%%%%%%%%%%%%%%%%%%%%%%%%%%%%%%%%%%%%%%%%%%%%%%%%%%%%%%%%%%%%%%%%%%
\section{Proof of Theorems \ref{thm:fatorodd-first} and \ref{thm:fatorodd-second} for $m\geq2$}
For all non-negative integers $a_1,\ldots,a_l$, and $k$, let
$$
C(a_1,\ldots,a_l;k)=\prod_{i=1}^l {a_i+a_{i+1}+1\brack a_i-k},
$$
where $a_{l+1}=a_1$, and let
\begin{align}
&\hskip -2mm S_r(n_1,\ldots,n_{m};j,q)  \notag \\[5pt]
&=\frac{(q;q)_{n_1}(q;q)_{n_m}}{(q;q)_{n_1+n_m+1}}
\sum_{k=0}^{n_1} q^{j(k^2+k)-(r+1)k}2k+1][k]^r[k+1]^r C(n_1,\ldots,n_{m};k),  \label{eq:sr-n1nm}\\[5pt]
&\hskip -2mm T_r(n_1,\ldots,n_{m};j,q)  \notag \\[5pt]
&=\frac{(q;q)_{n_1}(q;q)_{n_m}}{(q;q)_{n_1+n_m+1}}
\sum_{k=0}^{n_1}(-1)^k q^{{k\choose 2}+j(k^2+k)-rk} [2k+1][k]^r[k+1]^r C(n_1,\ldots,n_{m};k).  \label{eq:tr-n1nm}
\end{align}
It is easy to see that, for $m\geq 3$,
\begin{align}
C(n_1,\ldots,n_m;k)=\frac{(q;q)_{n_2+n_3+1}(q;q)_{n_m+n_1+1}}{(q;q)_{n_1+k+1}(q;q)_{n_2-k}(q;q)_{n_m+n_3+1}}
{n_1+n_2+1\brack n_1-k}C(n_3,\ldots,n_m;k).  \label{eq:C}
\end{align}
Applying \eqref{eq:C} and the $q$-Chu-Vandermonde identity (see, for example, \cite[p.~37, (3.3.10)]{Andrews})
\begin{align}
{n_1+n_2+1\brack n_1-k}
=\sum_{s=0}^{n_1-k}\frac{q^{s(s+2k+1)}(q;q)_{n_1+k+1}(q;q)_{n_2-k}}{(q;q)_s (q;q)_{s+2k+1}(q;q)_{n_1-k-s}(q;q)_{n_2-k-s}},  \label{eq:vandermonde}
\end{align}
we may write \eqref{eq:sr-n1nm} as
\begin{align*}
&\hskip -2mm S_r(n_1,\ldots,n_{m};j,q) \\
&=\frac{(q;q)_{n_2+n_3+1}(q;q)_{n_1} (q;q)_{n_m}}{(q;q)_{n_m+n_3+1}}\sum_{k=0}^{n_1}\sum_{s=0}^{n_1-k}
\frac{q^{j(k^2+k)-(r+1)k}[2k+1][k]^r[k+1]^r C(n_3,\ldots,n_{m};k)}{(q;q)_{s}(q;q)_{s+2k+1}(q;q)_{n_1-k-s}(q;q)_{n_2-k-s}}  \\
&=\frac{(q;q)_{n_2+n_3+1}(q;q)_{n_1} (q;q)_{n_m}}{(q;q)_{n_m+n_3+1}}\sum_{l=0}^{n_1} q^{l^2+l}\sum_{k=0}^{l}
\frac{q^{(j-1)(k^2+k)-(r+1)k}[2k+1][k]^r[k+1]^r C(n_3,\ldots,n_{m};k)}{(q;q)_{l-k}(q;q)_{l+k}(q;q)_{n_1-l}(q;q)_{n_2-l}},
\end{align*}
where $l=s+k$. Noticing that
$$
\frac{C(n_3,\ldots, n_{m};k)}{(q;q)_{l-k}(q;q)_{l+k+1}}
=\frac{(q;q)_{n_{m}+n_3+1}}{(q;q)_{n_3+l+1}(q;q)_{n_{m}+l+1}}C(l,n_3,\ldots, n_{m};k),
$$
we obtain
\begin{align}
S_r(n_1,\ldots,n_{m};j,q)
=\sum_{l=0}^{n_1} q^{l^2+l}{n_1\brack l}{n_2+n_3+1\brack n_2-l} S_r(l,n_3,\ldots,n_{m};j-1,q),\ m\geqslant 3. \label{eq:S-recsr}
\end{align}
Moreover, for $m=2$, applying \eqref{eq:vandermonde} we conclude
\begin{align}\label{eq:Sn1n2}
S_r(n_1,n_{2};j,q)=\sum_{l=0}^{n_1} q^{l^2+l}{n_1\brack l}{n_2\brack l}S_r(l;j-1,q).
\end{align}

Similarly, we have the following recurrence for \eqref{eq:tr-n1nm}:
\begin{align}
T_r(n_1,\ldots,n_{m};j,q)
&=\sum_{l=0}^{n_1} q^{l^2+l}{n_1\brack l}{n_2+n_3+1\brack n_2-l} T_r(l,n_3,\ldots,n_{m};j-1,q),\ m\geqslant 3,\\[5pt]
T_r(n_1,n_{2};j,q)
&=\sum_{l=0}^{n_1} q^{l^2+l}{n_1\brack l}{n_2\brack l}T_r(l;j-1,q).  \label{eq:Tn1n2}
\end{align}

We now proceed by induction on $m$. In section 4, we have proved that Theorem \ref{thm:fatorodd-second} holds for $m=1$.
Suppose that Theorem \ref{thm:fatorodd-second} is true for $m-1$ ($m\geq2$) and $0\leqslant j\leqslant m-1$.
By the induction hypothesis and the relation $[l]{n_1\brack l}=[n_1]{n_1-1\brack l-1}$, it is easy to check that
\begin{align*}
&{n_1\brack l}S_r(l,n_3,\ldots,n_{m};j,q)\equiv0 \mod [n_1]^{\min\{1,r\}}[n_m]^{\min\{1,{r\choose 2}\}},\\[5pt]
&{n_1\brack l}T_r(l,n_3,\ldots,n_{m};j,q)\equiv0 \mod [n_1]^{\min\{1,r\}}[n_m]^{\min\{1,r\}}
\end{align*}
for any non-negative integer $l$.
It follows from
\eqref{eq:S-recsr}--\eqref{eq:Tn1n2} that Theorem \ref{thm:fatorodd-second} holds for $m$ and $1\leqslant j\leqslant m$.
Applying the identity ${\alpha\brack k}_{q^{-1}}={\alpha\brack k}_{q}q^{k^2-\alpha k}$, we have
\begin{align*}
&S_r(n_1,\ldots,n_{m};0,q)=S_r(n_1,\ldots,n_{m};m,q^{-1}) q^{n_2+\cdots+n_{m-1}+n_1n_2+\cdots+n_{m-1}n_m-r},\\[5pt]
&T_r(n_1,\ldots,n_{m};0,q)=T_r(n_1,\ldots,n_{m};m-1,q^{-1}) q^{n_2+\cdots+n_{m-1}+n_1n_2+\cdots+n_{m-1}n_m-r}.
\end{align*}
Therefore, Theorem \ref{thm:fatorodd-second} also holds for $m$ and $j=0$. This completes the proof of Theorem \ref{thm:fatorodd-second}.
Similarly, we can prove Theorem \ref{thm:fatorodd-first} for $m\geqslant 2$.

\medskip
\noindent{\it Remark.} If we apply the following form of the $q$-Chu-Vandermonde identity
\begin{align*}
{n_1+n_2+1\brack n_1-k}
=\sum_{s=0}^{n_1-k}\frac{q^{(n_1-k-s)(n_2-k-s)}(q;q)_{n_1+k+1}(q;q)_{n_2-k}}{(q;q)_s (q;q)_{s+2k+1}(q;q)_{n_1-k-s}(q;q)_{n_2-k-s}},
\end{align*}
we have
\begin{align*}
S_r(n_1,\ldots,n_{m};j,q)
&=\sum_{l=0}^{n_1} q^{(n_1-l)(n_2-l)}{n_1\brack l}{n_2+n_3+1\brack n_2-l} S_r(l,n_3,\ldots,n_{m};j,q),\ m\geqslant 3,
\end{align*}
and so on.

%%%%%%%%%%%%%%%%%%%%%%%%%%%%%%%%%%%%%%%%%%%%%%%%%%%%%%%%%%%%%%%%%%%%%%%%%%%%%%%%%%%%%%%%%%%%%%%%%%%%%%%%%%%%%%%%%%%%%
\section{Proof of Theorem \ref{thm:q-ballot}}
Let $\Phi_n(q)$ be the $n$-th {\it cyclotomic polynomial} in $q$, i.e.,
\begin{align*}
\Phi_n(q):=\prod_{\substack{1\leqslant k\leqslant n\\ \gcd(n,k)=1}}(q-\zeta^k),
\end{align*}
where $\zeta$ is a $n$-th primitive root of unity. Let $\lfloor x\rfloor$ denote the greatest integer not exceeding $x$.
We will need the following result (see, for example, \cite[(10)]{KW} or \cite{CH,GZ06}).
\begin{prop}\label{prop:factor}
The $q$-binomial coefficient ${m\brack k}$ can be written as
$$
{m\brack k}=\prod_{d}\Phi_d(q),
$$
where $d$ ranges over all positive integers such that
$\lfloor k/d\rfloor+\lfloor (m-k)/d\rfloor<\lfloor m/d\rfloor$.
\end{prop}

We now suppose that $r+s\equiv 1\pmod 2$ and $0\leqslant j\leqslant s$. Letting $m=s$ and $n_1=\cdots=n_s=n$ in \eqref{eq:first-1}, one sees that
\begin{align*}
&\sum_{k=0}^{n}q^{j(k^2+k)-(r+s)k} [2k+1]^{r+s}{2n+1\brack n-k}^s\equiv 0\mod [2n+1]{2n\brack n}.
%&\sum_{k=0}^{n}(-1)^k q^{{k\choose 2}+j(k^2+k)-(r+s-1)k} [2k+1]^{r+s}{2n+1\brack n-k}^s  \equiv 0 \mod [2n+1]{2n\brack n}.
\end{align*}
Noticing that
\begin{align}
[2k+1]{2n+1\brack n-k}q^{n-k}=[2n+1]\left({2n\brack n-k}-{2n\brack n-k-1}\right)\equiv 0 \mod [2n+1],  \label{eq:div-qballot}
\end{align}
we immediately get
\begin{align*}
&\sum_{k=0}^{n}q^{j(k^2+k)-rk}[2k+1]^r\left({2n\brack n-k}-{2n\brack n-k-1}\right)^s \equiv 0\mod \frac{{2n\brack n}}{\gcd\left({2n\brack n},[2n+1]^{s-1}\right)}.
%&\sum_{k=0}^{n} (-1)^kq^{{k\choose 2}+j(k^2+k)-(r-1)k}[2k+1]^r\left({2n\brack n-k}-{2n\brack n-k-1}\right)^s \frac{{2n\brack n}}{\gcd\left({2n\brack n},[2n+1]^{s-1}\right)}.
\end{align*}
But, by Proposition \ref{prop:factor} we have
\begin{align}
\gcd\left({2n\brack n},[2n+1]\right)=1.  \label{eq:gcd-bino}
\end{align}
This completes the proof of \eqref{eq:q-ballt-1}. Similarly, we can prove \eqref{eq:q-ballt-2}.

\medskip
\noindent{\it Remark.} In general, for any positive integer $n$, we cannot expect $\gcd({2n\choose n},2n+1)=1$.
This means that sometimes the $q$-analogue of a mathematical problem will be easier than the original one,
although in most cases the former will be much more difficult.

%%%%%%%%%%%%%%%%%%%%%%%%%%%%%%%%%%%%%%%%%%%%%%%%%%%%%%%%%%%%%%%%%%%%%%%%%%%%%%%%%%%%%%%%%%%%%%%%%%%%%%%%%%%%%%%%%%%%%
\section{Proof of Theorem \ref{thm:mn-qq}}
We first give the following result, which is a generalization of \eqref{eq:gcd-bino}.
\begin{lem}\label{lem:gcd-2m2n}
For all $m,n\in\Z$, there holds
\begin{align}
\gcd\left(\frac{[2m]![2n]!}{[m+n]![m]![n]!},[2m+1]\right)=1.  \label{eq:gcd-2m2n}
\end{align}
\end{lem}
\pf It is well known that
$$
q^n-1=\prod_{d|n}\Phi_d(q),
$$
and so
$$
[n]!=(q-1)^{-n}\prod_{k=1}^{n}(q^k-1)=(q-1)^{-n}\prod_{d=1}^{n}\Phi_d(q)^{\lfloor\frac{n}{d}\rfloor}.
$$
Therefore,
\begin{align*}
\frac{[2m]![2n]!}{[m+n]![m]![n]!}
=\prod_{d=1}^{\max\{2m,2n\}}\Phi_d(q)^{\lfloor\frac{2m}{d}\rfloor+\lfloor\frac{2n}{d}\rfloor-\lfloor\frac{m+n}{d}\rfloor-\lfloor\frac{m}{d}\rfloor-\lfloor\frac{n}{d}\rfloor}.
\end{align*}

For any irreducible factor $\Phi_d(q)$ of $[2m+1]$, we have $2m+1\equiv 0 \pmod d$. It follows that $d$ is odd and $m\equiv \frac{d-1}{2}\pmod{d}$.
Suppose that $n\equiv a\pmod{d}$ with $0\leqslant a\leqslant d-1$. We consider the following two cases. If $a\leqslant \frac{d-1}{2}$, then
\begin{align}
&\left\lfloor\frac{2m}{d}\right\rfloor+\left\lfloor\frac{2n}{d}\right\rfloor
-\left\lfloor\frac{m+n}{d}\right\rfloor-\left\lfloor\frac{m}{d}\right\rfloor-\left\lfloor\frac{n}{d}\right\rfloor  \notag\\[5pt]
&\quad=\frac{2m-d+1}{d}+\frac{2n-2a}{d}-\frac{m+n-\frac{d-1}{2}-a}{d}-\frac{m-\frac{d-1}{2}}{d}-\frac{n-a}{d} \notag\\[5pt]
&\quad=0.  \label{eq:fraction}
\end{align}
If $a\geqslant \frac{d+1}{2}$, then the left-hand side of \eqref{eq:fraction} is equal to
$$
\frac{2m-d+1}{d}+\frac{2n-2a+d}{d}-\frac{m+n+\frac{d+1}{2}-a}{d}-\frac{m-\frac{d-1}{2}}{d}-\frac{n-a}{d}=0.
$$
This means that $\Phi_d(q)$ is not a factor of $\frac{[2m]![2n]!}{[m+n]![m]![n]!}$, and so the formula \eqref{eq:gcd-2m2n} holds.
\qed

It is clear that Theorem~\ref{thm:fatorodd-first} can be restated as follows.

\begin{thm}\label{thm:req-odd}
Let $n_1,\ldots,n_m\in\Z$ and $j,r\in\N$ with $j\leqslant m$. Then the expressions
\begin{align}
[n_1]!\prod_{i=1}^m\frac{[n_i+n_{i+1}+1]!}{[2n_i+1]!}
\sum_{k=0}^{n_1} q^{j(k^2+k)-(2r+1)k} [2k+1]^{2r+1}\prod_{i=1}^m {2n_i+1\brack n_i-k}, \\[5pt]
[n_1]!\prod_{i=1}^m\frac{[n_i+n_{i+1}+1]!}{[2n_i+1]!}
\sum_{k=0}^{n_1} (-1)^k q^{{k\choose 2}+j(k^2+k)-2rk} [2k+1]^{2r+1}\prod_{i=1}^m {2n_i+1\brack n_i-k}
\end{align}
where $n_{m+1}=-1$, are Laurent polynomials in $q$ with integer coefficients.
\end{thm}

\noindent{\it Proof of Theorem \ref{thm:mn-qq}.}
Letting $n_1=\cdots= n_s=m$ and
$n_{s+1}=\cdots=n_{s+t}=n$ in Theorem~\ref{thm:fatorodd-first}, we obtain
\begin{align}
&[m+n+1]\sum_{k=0}^m q^{j(k^2+k)-(r+s+t)k} [2k+1]^{r+s+t}{2m+1\brack m-k}^s{2n+1\brack n-k}^t  \notag\\[5pt]
&\quad \equiv 0 \mod \frac{[2m+1]![2n+1]!}{[m+n]![m]![n]!}.  \label{eq:big-sum-mn}
\end{align}
By \eqref{eq:div-qballot} and the definition of $q$-ballot numbers $A_{n,k}(q)$, we deduce from \eqref{eq:big-sum-mn} that
\begin{align*}
&[m+n+1]\sum_{k=0}^m q^{j(k^2+k)-rk} [2k+1]^{r}A_{m,k}(q)^s  A_{n,k}(q)^t  \\[5pt]
&\quad\equiv 0 \mod \frac{\frac{[2m]![2n]!}{[m+n]![m]![n]!}}{\gcd(\frac{[2m]![2n]!}{[m+n]![m]![n]!},[2m+1]^{s-1}[2n+1]^{t-1})}.
\end{align*}
By Lemma \ref{lem:gcd-2m2n}, we have
$$
\gcd\left(\frac{[2m]![2n]!}{[m+n]![m]![n]!},[2m+1]^{s-1}[2n+1]^{t-1}\right)=1.
$$
This completes the proof.  \qed

Letting $m=n+1$ or $m=2n$ in Theorem \ref{thm:mn-qq}, we get the following result, which in the $q=1$ case confirms
a conjecture of Guo and Zeng \cite[Conjecture 6.10]{GZ2011}. Note that $\frac{1}{[n+1]}{2n\brack n}$ is the famous
{\it $q$-Catalan numbers} (see \cite{FH}).
\begin{cor}\label{cor:n+4and4n}
Let $n,s,t\in\Z$ and $j,r\in \N$ with $r+s+t\equiv 1\pmod 2$ and $j\leqslant s+t$. Then
\begin{align*}
\sum_{k=0}^n \tau_k [2k+1]^{r}A_{n+1,k}(q)^s  A_{n,k}(q)^t
&\equiv 0 \mod \frac{1}{[n+1]}{2n\brack n}, \\[5pt]
\sum_{k=0}^n \tau_k [2k+1]^{r}A_{2n,k}(q)^s  A_{n,k}(q)^t
&\equiv 0 \mod \frac{1}{[3n+1]}{4n\brack n},
\end{align*}
where $\tau_k=q^{j(k^2+k)-rk}$ or $\tau_k=(-1)^kq^{{k\choose 2}+j(k^2+k)-(r-1)k}$.
\end{cor}

%%%%%%%%%%%%%%%%%%%%%%%%%%%%%%%%%%%%%%%%%%%%%%%%%%%%%%%%%%%%%%%%%%%%%%%%%%%%%%%%%%%%%%%%%%%%%%%%%%%%%%%%%%%%%%%%%%%%%
\section{Some consequences and conjectures}
In this section, we will give some consequences of Theorem \ref{thm:fatorodd-first}. Most of these results are $q$-analogues of
the corresponding results listed in \cite[Section 6]{GZ2011}.
 Note that there are exactly similar consequences of Theorem \ref{thm:fatorodd-second}.
We shall also confirm some conjectures in \cite[Section 6]{GZ2011}. For convenience,
we let $\varepsilon_k=q^{j(k^2+k)-(2r+1)k}$ or $\varepsilon_k=(-1)^kq^{{k\choose 2}+j(k^2+k)-2rk}$ throughout this section.

Letting $n_{2i-1}=m$ and $n_{2i}=n$ for $i=1,\ldots,a$ in Theorem \ref{thm:fatorodd-first} and observing the symmetry of $m$ and $n$,
we obtain
\begin{cor}\label{cor:mnrs}
Let $a,m,n\in\Z$ and $j,r\in\N$ with $j\leqslant 2a$. Then
\begin{align*}
\sum_{k=0}^{m}\varepsilon_k[2k+1]^{2r+1} {m+n+1\brack m-k}^a {m+n+1\brack n-k}^a
&\equiv 0\mod [m+n+1]{m+n\brack m}.
\end{align*}
\end{cor}

Letting $n_{3i-2}=l$, $n_{3i-1}=m$ and $n_{3i}=n$ for $i=1,\ldots,a$ in Theorem~\ref{thm:fatorodd-first}, we get
\begin{cor}
Let $a,l,m,n\in\Z$ and $j,r\in\N$ with $j\leqslant 3a$. Then
\begin{align*}
\sum_{k=0}^{m}\varepsilon_k[2k+1]^{2r+1} {l+m+1\brack l-k}^a {m+n+1\brack m-k}^a {n+l+1\brack n-k}^a
\equiv 0\mod [m+n+1]{m+n\brack m}.
\end{align*}
\end{cor}

Taking $m=2a+b$ and letting $n_i=n$ if $i=1,3,\ldots,2a-1$ and $n_i=n-1$ otherwise in
Theorem~\ref{thm:fatorodd-first}, we get
\begin{cor}
Let $a,n\in\Z$ and $b,j,r\in\N$ with $j\leqslant 2a+b$. Then
\begin{align*}
\sum_{k=0}^{n-1}\varepsilon_k[2k+1]^{2r+1} {2n\brack n-k}^a {2n\brack n-k-1}^a {2n-1\brack n-k-1}^b
\equiv 0\mod [n]{2n\brack n}.
\end{align*}
\end{cor}

By Theorem~\ref{thm:req-odd} it is easily seen that, for all $a_1,\ldots,a_m\in\Z$,
\begin{align}
[n_1]!\prod_{i=1}^m\frac{[n_i+n_{i+1}+1]!}{[2n_i+1]!}
\sum_{k=0}^{n_1} \varepsilon_k [2k+1]^{2r+1}\prod_{i=1}^m {2n_i+1\brack n_i-k}^{a_i}\quad(n_{m+1}=-1) \label{eq:req-power}
\end{align}
is a Laurent polynomial in $q$ with integer coefficients. For $m=3$, letting $(n_1,n_2,n_3)$ be $(n,n+2,n+1)$, $(n,3n,2n)$, $(2n,n,3n)$,
$(2n,n,4n)$, or $(3n,2n,4n)$, we immediately get the following three conclusions.
\begin{cor}\label{cor:n1n2n3}
Let $a,b,c,n\in\Z$ and $j,r\in \N$ with $j\leqslant a+b+c$. Then
\begin{align}
\sum_{k=0}^n \varepsilon_k [2k+1]^{2r+1}   {2n+1\brack n-k}^a {2n+3\brack n-k+1}^b {2n+5\brack n-k+2}^c
&\equiv 0 \mod [2n+5]{2n+1\brack n}. \label{eq:n1n2n3}
\end{align}
\end{cor}
\begin{cor}\label{cor:qfactor-246n}
Let $a,b,c,n\in\Z$ and $j,r\in \N$ with $j\leqslant a+b+c$. Then
\begin{align*}
\sum_{k=0}^n \varepsilon_k [2k+1]^{2r+1} {6n+1\brack 3n-k}^a {4n+1\brack 2n-k}^b {2n+1\brack n-k}^c
&\equiv 0 \mod [2n+1]{6n+1\brack n}, \\
\sum_{k=0}^n \varepsilon_k [2k+1]^{2r+1} {6n+1\brack 3n-k}^a {4n+1\brack 2n-k}^b {2n+1\brack n-k}^c
&\equiv 0 \mod [2n+1]{6n+1\brack 3n}.
\end{align*}
\end{cor}

\begin{cor}\label{cor:qfactor-248n}
Let $a,b,c,n\in\Z$ and $j,r\in \N$ with $j\leqslant a+b+c$. Then
\begin{align*}
[3n+1]\sum_{k=0}^n \varepsilon_k [2k+1]^{2r+1} {8n+1\brack 4n-k}^a {4n+1\brack 2n-k}^b {2n+1\brack n-k}^c
&\equiv 0 \mod [2n+1][4n+1]{8n+1\brack 3n},   \\
\sum_{k=0}^n  \varepsilon_k [2k+1]^{2r+1} {8n+1\brack 4n-k}^a {6n+1\brack 3n-k}^b {4n+1\brack 2n-k}^c
&\equiv 0 \mod [4n+1]{8n+1\brack 3n},
\end{align*}
\end{cor}

We have the following conjectural generalization of Corollaries \ref{cor:qfactor-246n} and \ref{cor:qfactor-248n}.
\begin{conj}
Let $n,r,s,t\in\Z$ with $r+s+t\equiv 1\pmod 2$ and $j\in\N$. Then
\begin{align*}
[4n+1]\sum_{k=0}^n \eta_k A_{3n,k}(q)^r A_{2n,k}(q)^s A_{n,k}(q)^t
&\equiv 0 \mod \frac{1}{[6n+1]}{6n+1\brack n}, \\
[4n+1]\sum_{k=0}^n \eta_k A_{3n,k}(q)^r A_{2n,k}(q)^s A_{n,k}(q)^t
&\equiv 0 \mod \frac{1}{[6n+1]}{6n+1\brack 3n}, \\
[8n+1]\sum_{k=0}^n \eta_k A_{4n,k}(q)^r A_{2n,k}(q)^s A_{n,k}(q)^t
&\equiv 0 \mod{ {8n+1\brack 3n}}, \\
[6n+1][8n+1]\sum_{k=0}^n \eta_k A_{4n,k}(q)^r A_{3n,k}(q)^s A_{2n,k}(q)^t
&\equiv 0 \mod{ {8n+1\brack 3n} },
\end{align*}
where $\eta_k=q^{j(k^2+k)}$ or $\eta_k=(-1)^kq^{{k+1\choose 2}+j(k^2+k)}$.
\end{conj}

For general $m\geq 2$, in \eqref{eq:req-power} taking $(n_1,\ldots,n_{m})$ to be
\begin{align*}
\begin{cases}
(n,n+2,\ldots,n+m-1,n+m-2,n+m-4,\ldots,n+1),&\text{if $m$ is odd,}\\[5pt]
(n+1,n+3,\ldots,n+m-1,n+m-2,n+m-4,\ldots,n),&\text{if $m$ is even,}
\end{cases}
\end{align*}
we are led to the following generalization of \eqref{eq:n1n2n3}.
\begin{cor}\label{cor:final}
Let $m\geq 2$, and let $n,a_1,\ldots,a_m\in\Z$ and $j,r\in\N$ with $j\leqslant a_1+\cdots+a_m$. Then
\begin{align*}
\sum_{k=0}^n \varepsilon_k [2k+1]^{2r+1} \prod_{i=1}^m {2n+2i-1\brack n+i-k-1}^{a_i}
&\equiv 0 \mod [2n+2m-1]{2n+1\brack n}.
\end{align*}
\end{cor}
We have the following challenging conjecture related to Corollary \ref{cor:final}.
\begin{conj}\label{conj:final}
Let $n,r_1,\ldots,r_m\in\Z$ with $r_1+\cdots+r_m\equiv 1\pmod 2$ and $j\in\N$, there holds
\begin{align*}
\sum_{k=0}^n \eta_k \prod_{i=1}^m A_{n+i-1,k}(q)^{r_i}
&\equiv 0 \mod \frac{1}{[n+1]}{2n\brack n},
\end{align*}
where $\eta_k=q^{j(k^2+k)}$ or $\eta_k=(-1)^kq^{{k+1\choose 2}+j(k^2+k)}$.
\end{conj}
Note that, for $m=1$ and $0\leqslant j\leqslant r_1$, Conjecture \ref{conj:final} is true
by Theorem \ref{thm:q-ballot}. For $m=2$ and $0\leqslant j\leqslant r_1+r_2$,
Conjecture \ref{conj:final} is also true by the first congruence in Corollary~\ref{cor:n+4and4n}.
Note that the $q=1$ case of Conjecture \ref{conj:final} has been checked by Guo and Zeng \cite{GZ2011} for $n=2$, or $m\leq 6$ and
$n=4,9,10,11,3280,7651,7652$.

We end the paper with the following conjecture.
\begin{conj}Theorems {\rm\ref{thm:fatorodd-first}} and {\rm\ref{thm:fatorodd-second}} hold for all $j\in\N$.
\end{conj}

\vskip 5mm \noindent{\bf Acknowledgments.} The first author was partially
supported by the National Natural Science Foundation of China (grant 11371144),
the Natural Science Foundation of Jiangsu Province (grant BK20161304),
and the Qing Lan Project of Education Committee of Jiangsu Province.

\renewcommand{\baselinestretch}{1}

\end{document}